\input amstex
\documentstyle{amsppt}
\magnification1200
\tolerance=10000
\def\n#1{\Bbb #1}

\def\e11{E_{11}}

\def\ga{\goth A}

\def\ga{\gamma}
\def\la{\lambda}

\def\be{\beta}

\def\g{\goth }
\def\al{\alpha}

\def\th{\theta}
\def\vk{\varkappa}

\NoRunningHeads
\topmatter
\title
Analytic ranks of twists of Carlitz modules --- a survey of results
\endtitle
\author
A. Grishkov, D. Logachev\footnotemark \footnotetext{E-mails: shuragri{\@}gmail.com; logachev94{\@}gmail.com (corresponding author)\phantom{****************}}
\endauthor
\address
First author: Departamento de Matem\'atica e estatistica
Universidade de S\~ao Paulo. Rua de Mat\~ao 1010, CEP 05508-090, S\~ao Paulo, Brasil, and Omsk State University n.a. F.M.Dostoevskii. Pr. Mira 55-A, Omsk 644077, Russia.
\medskip
Second author: Departamento de Matem\'atica, Universidade Federal do Amazonas, Manaus, Brasil
\endaddress
\abstract We give in this paper a survey of results obtained in our earlier papers, and state explicitly some problems of further research, for example: are the analytic ranks bounded, or not? Twists of Carlitz modules are parametrized by polynomials over finite fields $\n F_q$. The analytic rank of a twist is the order of zero of its L-function at a point. The set of polynomials of degree $\le m$ such that the analytic ranks of the corresponding twists are $\ge i$ is $X(m,i)(\n F_q)$ where  $X(m,i)$ is an affine variety defined over $\n F_p$ (we do not know what is its dimension). We consider also a related invariant of a twist, namely, the behaviour of its L-function at infinity (the rank at infinity). We know much more on varieties corresponding to twists of a fixed rank at infinity and on their lifts from $\n F_p$ to $\n Z$ . For example, for $q=2$ the irreducible components of these varieties are described in terms of finite rooted weighted binary trees. A similar description for $q>2$ is not found yet. 
\endabstract
\endtopmatter
\document
{\bf 1. Introduction.}
\medskip
The authors started to study in [GL16], [GLZ22], [GLZ23] and [GLZ] the behaviour of the analytic rank ( = orders of vanishing of $L$-functions) of twists of Carlitz modules (the simplest analogs of elliptic curves in finite characteristic). Also, the authors studied in these papers some algebraic varieties over both $\n F_p$ and $\n Q$ related to the sets of twists of a fixed rank. There are many open problems in this theory, hence many possibilities of further research. Unfortunately all the above papers are long and overwhelmed with technical details that makes their reading a tedious task. 
\medskip
Hence, in the present paper we give a concise description of the results of the above papers. This description will permit to the reader to understand what was made, and to continue this research. We think that this research will be interesting, because it will give relations of some different branches of mathematics: the theory of determinantal varieties, graph theory and combinatorics (irreducible components of some attached varieties are described in terms of finite rooted weighted binary trees and their unions; this result was obtained by a kind of "inspiration", and not by a routine development of the theory). The research problems are stated in Problems 1.4; 3.4; Section 3.6; Conjecture 4.1; Problem 4.6; Section 6. 
\medskip
It turns out that the main formula (3.2) and the Theorem 3.3 giving an explicit expression of a $L$-function under consideration, reduces the whole theory to the problems of algebraic geometry. It is sufficient for the reader to use the formula (3.2). For the development of the theory he should not know what is an Anderson t-motive, a Carlitz module, and what is its twist - the further research is reduced enturely to methods of algebraic geometry. 
\medskip
Twists of Carlitz modules are Anderson t-motives. There are several types of $L$-functions of a t-motive $M$; we consider the simplest of them. It is denoted by $L(M,T)$. 
\medskip
Later we consider only the case when $M$ is a twist of a Carlitz module. For these $M$ we have:
$$L(M,T)\in(\n F_q[t])[T]\eqno{(1.1)}$$
where $T$ is an independent variable (it is an analog of
$p^{-s}$ of the classical $L$-function),  $\n F_q[t]$ is the ring of coefficients of the ring of polynomials in $T$; here $t$ is also an independent variable. 
\medskip
{\bf Remark 1.2.} For any Anderson t-motive $M$ (non necessary a twist of a Carlitz module) we have $L(M,T)\in(\n F_q[t])(T)$. Hence, the analytic rank of $M$ (see the below definition) is defined for all t-motives $M$.
\medskip
We define the analytic rank of $M$, by analogy with the number field case, as the order of zero of $L(M,T)$ at some value of $T$. Unlike the number field case, there is no functional equation for $L(M,T)$. What value of $T$ should we choose for the definition of the analytic rank? 
\medskip
For $T=0$ we have always $L(M,T)=1$, hence we choose the value $T=1$ (clearly we consider 1 as an element of $\n F_q[t]$; the case of $T=a$ for any $a\in \n F_q$ is immediately reduced to the case $a=1$, see [GL16], Remark 4.3).
\medskip
Hence, we state 
\medskip
{\bf Definition 1.3.} The analytic rank of $M$ is the order of vanishing of $L(M,T)$ at $T=1$. It is denoted by $r(M)$. 
\medskip
Can the reader offer a better definition of the analytic rank?
\medskip
Since $L(M,T)$ is a polynomial in $T$ with coefficients in $\n F_q[t]$, the condition that the analytic rank of $M$ is a fixed number $r$ is equivalent to the conditions that some elements of $\n F_q[t]$ (coming from the coefficients of $L(M,T)$ are identically zero. 
\medskip
As we mentioned above, it is not necessary for the reader to know what a Carlitz module, and its twists, are. Let us give the only information on these objects that he should know. 
\medskip
Let $q$ be a power of a prime $p$. There exists the $q$-Carlitz module denoted by $\g C=\g C(q)$. The set of twists of $\g C$ is an abelian group $(\n F_q(\th))^*/(\n F_q(\th))^{*(q-1)}$ where $\th$ is a free variable; we denote it by $\g T$. The representatives of cosets of $(\n F_q(\th))^{*(q-1)}$ in $(\n F_q(\th))^*$ ( = elements of $\g T$) are polynomials $P$ free from $(q-1)$-th powers. The twist of $\g C$ corresponding to such $P$ is denoted by $\g C_P$, and its rank is denoted by $r(\g C_P)$. 
\medskip
We see that for $q=2$ all twists are isomorphic, hence the simplest non-trivial case is $q=3$. If $P_1, \ P_2\in \n F_q(\th)^*$ and $P_1/P_2 \in (\n F_q(\th))^{*(q-1)}$ then the corresponding twists $\g C_{P_1}$, $\g C_{P_2}$ are isomorphic, they have equal ranks, and their $L$-functions differ by finitely many Euler factors corresponding to the primes of bad reduction. See [GL16], (5.6.1) for an explicit formula. 
\medskip
Since Drinfeld modules are analogs of elliptic curves, we shall recall the corresponding facts from the theory of elliptic curves. Namely, let $E$ be an elliptic curve over $\n Q$. Its twists $E_d$ are parametrized by squarefree $d\in \n Z$. The $(q-1)$-power-free $P\in \n F_q[\th]-\{0\}$ are the function field case analogs of squarefree $d\in \n Z$.
\medskip
In addition to study of the behaviour of $L(\goth C_P,T)$ at $T=1$, we are interested to study the behaviour of $L(\goth C_P,T)$ at $T=\infty$, i.e. the degree of $L(\goth C_P,T)$ as a polynomial in $T$. We start to study this invariant in [GL16], Part II; the papers [GLZ22], [GLZ23], [GLZ] are entirely devoted to its study. 
\medskip
The principal open problem is the following: 
\medskip
{\bf Problem 1.4.} Let $q$ be fixed. Are $r(\g C_P)$ bounded? If yes, what is 
\medskip
(a) The maximal value of $r(\g C_P)$?
\medskip
(b) The maximal value of $r(\g C_P)$ which occurs for infinitely many $P$ free from $(q-1)$-th powers? 
\medskip
A related important open problem is the Problem 3.4, see below.
\medskip
{\bf Example 1.5.} For $q=3$ the naive considerations of type "quantity of variables minus quantity of equations" (see 3.5 below for the details) predict that the maximal value of $r(\g C_P)$ is 3. Really, there exist two squarefree polynomials $P$ of degree 39 having $r(\g C_P)=6$ (see [GL16], Table 6.11). 
\medskip
Since $L(\goth C_P,T)$ "consists of polynomials", the sets of $P$ having the analytic rank $\ge i$, where $i\ge0$ is fixed, are algebraic varieties (more exactly, $\n F_q$-points of some ind-varieties over $\n F_q$). We shall also consider the natural lifts of these polynomials in characteristic 0, and we shall study the corresponding algebraic varieties over $\n Q$. 
\medskip
{\bf 2. More information.} 
\medskip
Let us recall the behaviour of ranks of twists of an elliptic curve $E$ over $\n Q$. The parity of the rank of $E_d$ is the sign of its functional equation, it depends on a residue of $d$ by some module. For almost all $d$ the rank of $E_d$ takes the minimal possible value, i.e. it is 0 for even curves and 1 for the odd ones. Rare jumps of the rank of $E_d$ occur, i.e. occasionally the rank of an even $E_d$ can be 2 (rarely), 4 (much more rarely) etc, and respectively the rank of an odd $E_d$ can be 3 (rarely), 5 (much more rarely) etc. 
\medskip
Are the ranks bounded? [PPVW] suggests that yes. These jumps of rank seem to be irregular, i.e. we do not know a formula for $d$ such that the rank of $E_d$ has a jump. 
\medskip
The behavior of $r(\g C_P)$ is the following. The group of twists $\g T$ contains a subgroup $S_0$ of index $(q-1)^2$ (see [GL16], Proposition 4.5 and the above lines). There exists a coset $S$ of  $S_0$ such that for the majority of $P\not\in S$ we have $r(\g C_P)=0$, for the majority of $P\in S$ we have $r(\g C_P)=1$. Further, for both $P\in S, \ P\not\in S$ rare jumps of $r(\g C_P)$ occur. The parity of $r(\g C_P)$ can be arbitrary; also, the residues modulo $q-1$ of $r(\g C_P)$ can be arbitrary (recall that the residue modulo $q-1$ is the analog of parity for the function field case). 
\medskip
Let us give more information. Let $m$ be fixed. We consider the set of polynomials $P\in \n F_q[\th]$ of degree $\le m$, not necessarily free from $(q-1)$-th powers: 
$$P=\sum_{k=0}^m a_k\theta^k \hbox{ where } a_k\in \n F_q.\eqno{(2.1)}$$

Hence, we can identify the set of these $P$ with the set $\n A^{m+1}(\n F_q):=\n F_q^{m+1}$ --- the affine space $\n A^{m+1}(\n F_q)$ over $\n F_q$ of dimension $m+1$. 
\medskip
{\bf Theorem 2.2.} For any $i\ge0$ there exists a {\bf naturally defined} affine variety $$X(m,i)\subset \n A^{m+1}(\overline{\n F_q})$$ defined over $\n F_p$ such that the set of $P$ of degree $\le m$ having $r(\g C_P)\ge i$ is $X(m,i)(\n F_q)$.
\medskip
Further, for any $m$ we have $X(m,i)=X(m+1,i)\cap  \n A^{m+1}(\overline{\n F_q})$ in $\n A^{m+2}(\overline{\n F_q})$. Hence, $$\underset{m\ge1}\to{\bigcup}X(m,i)$$ is an ind-variety. 
\medskip
{\bf Remark 2.3.}\footnotemark \footnotetext{The authors are grateful to anonymous reviewer who indicated them this information.} Since the set of $P$ of degree $\le m$ having $r(\g C_P)\ge i$ is a finite set, it is always the set of 
$\n F_q$-points of an affine variety. Hence, the key words of Theorem 2.2 are: 
\medskip
$X(m,i)$ is a {\it naturally defined} variety. 
\medskip
Its definition comes from the proof of Theorem 3.3. 
\medskip
So, in future we shall consider exactly $X(m,i)$ from Theorem 2.2.
\medskip
Clearly $X(m,i)=X(q,m,i)$ depends on $q$.
\medskip
{\bf 3. Formula for $L(\g C_P,T)$.} 
\medskip
It is convenient to distinguish in the below formulas the abstact variables $\g a_1, \dots, \g a_m$ and their explicit values $a_0, \dots, a_m\in \n F_q$ obtained by assignments $\g a_i \mapsto a_i$ (unfortunately this was not made in [GL16], [GLZ22], [GLZ23], [GLZ]).  
\medskip
Let $q$ be fixed, and let $$k=\lceil\frac{m+1}{q-1}\rceil\eqno{(3.1)}$$where $\lceil x \rceil:=\min \{n\in \n Z\ |\ n\ge x\}$ is the ceiling function. There exists a $k\times k$-matrix denoted by $\g M(\g a_0, \dots, \g a_m)$. See [GL16], Formula (3.1), case $n=1$ for its definition, and [GL16], Formula (3.2) for its explicit form for the case $n=1$ --- the only case that we need now. It is a version of a Schur matrix\footnotemark \footnotetext{See, for example, [FP], (1.5) for a definition of a Schur matrix.}, but its entries belong to $\n Z[\g a_0, \dots, \g a_m,t]$, i.e. they depend not only on $\g a_0, \dots, \g a_m$, but on a variable $t$ as well. We consider 

$$L_m:=\det(I_k - \g M(\g a_0, \dots, \g a_m)T)\in \n Z[\g a_0, \dots, \g a_m,t,T].\eqno{(3.2)}$$
\medskip
{\bf Theorem 3.3} ( = [GL16], Theorem 3.3). For $P$ from (2.1) we have: $L(\g C_P,T)$ is the specialization of $L_m$ obtained by the assignments $\forall \ i$ $\g a_i \mapsto a_i$, and reduction of coefficients of $L_m$ from $\n Z$ to $\n F_p$. 
\medskip
For any $\al,\ \be$ we denote the coefficient of $L_m$ at $t^\al T^\be$ by $H_{\be\al}$ (unfortunately, there is no concordance of this notations in different sources). We have $H_{\be\al}\in \n Z[\g a_0, \dots, \g a_m]$; they depend also on $m$ and $q$. We denote their reduction in $\n F_p$ by $\overline{H_{\be\al}}$. 
\medskip
Theorem 3.3 implies immediately that $X(m,i)\subset \n A^{m+1}(\n F_q)$ is the set of zeroes of some polynomials that are linear combinations of $\overline{H_{\be\al}}$ with coefficients in $\n F_p$. 
\medskip
We have the following open problem (recall that $q$ is fixed):
\medskip
{\bf Problem 3.4.} What are the properties of $X(m,i)$ (singularities, irreducible components etc.)? What are their dimensions? What is the minimal value of $i$ such that $X(m,i)=\emptyset$? What is the asymptotics of $\#(X(m,i)(\n F_q))$ as $m$ grows?
\medskip
Clearly if $\forall \ m $ we have $X(m,i)=\emptyset$ then the maximal value of the rank of twists of $\g C$ is $< i$, but the converse is not necessarily true. 
\medskip
We know nothing about $X(m,i)$. The only information on their dimension is the following. 
We have $m+1$ variables $\g a_0, \dots, \g a_m$ and some quantity of the above linear combinations of $\overline{H_{\be\al}}$. 
\medskip
{\bf 3.5.} The naive formula for the dimension "$m+1$ minus the quantity of equations" (i.e. as $X(m,i)$ were a complete intersection) is given in [GL16], Proposition 6.6. For $q=3$ it gives that $X(m,i)=\emptyset$ for $i\ge4$. 
\medskip
To check, whether $X(m,i)$ are complete intersections, or not, we calculated $\#(X(m,i)(\n F_3))$ for $q=3$ and for some values of $m$ and $i$ (see tables in [GL16], Section 6). These values show that most likely $X(m,i)$ are not complete intersections. As we mentioned above, there are two squarefree $P$ of degree 39 having $r(\g C_P)=6$. 
\medskip
{\bf 3.6.} There are many questions related to these tables. 
\medskip
{\bf 3.6.1.} Why high values of $r(\g C_P)$ occur for shift-stable\footnotemark \footnotetext{$P\in \n F_q[\th]$ is called shift-stable if $\forall \ a\in \n F_q$ we have $P(\th+a)=P(\th)$, i.e. $P$ is a polynomial in $\th^q-\th$. For $q=3$ all known $P$ of ranks 5, 6 are shift-stable.} $P$?  
\medskip
{\bf 3.6.2.} Why for $P$ of ranks 5, 6 (four such polynomials were found) we have $L(\goth C_P,T)\in \n F_3[T]$, although a priori $L(\goth C_P,T)\in (\n F_q[t])[T]$? 
\medskip
{\bf 3.6.3.} Why for $P$ of ranks 5, 6 we have: $L(\goth C_P,T)$ has a factorization in factors of small degrees?
\medskip
{\bf 3.6.4.} Why there exists a correlation between high values of $r(\goth C_P)$ and $r_\infty(\goth C_P)$ (the rank at infinity, see Section 4)? 
\medskip
Moreover, it seems that high-rank polynomials occur among polynomials which are stable with respect to "sufficiently large" subgroups of $GL_2(\n F_q)$ (see [GL16], Section 2 for the action of $GL_2(\n F_q)$ on polynomials $P$). 
\medskip
{\bf 4. Rank at infinity.}
\medskip
The rank at infinity of $\g C_P$, denoted by $r_\infty(\g C_P)$, indicates the behaviour of $L(\goth C_P,T)$ at $T=\infty$. Namely, it is (here $k$ is from (3.1))
$$k-\deg_T(L(\g C_P,T))$$
i.e. the deficiency of the degree of $L(\g C_P,T)$ with respect to the maximal possible degree (according (3.2), the maximal possible degree of $L(\goth C_P,T)$ is $k$; all degrees are considered with respect to the variable $T$). 
\medskip
For any $m,\ i$ we denote by $X_\infty(m,i)$ the set of $P$ of degree $\le m$ (or, equivalently, the set of elements of $\n A^{m+1}$) such that $r_\infty(\g C_P)\ge i$. 
\medskip
We have: $X_\infty(m,i)$ is the set of zeroes of $\overline{H_{\be\al}}$, where $\be>k-i$ (this is immediate). 

Polynomials $H_{\be\al}$ are homogeneous, hence we can consider the projectivization of $X_\infty(m,i)$ (denoted by $X_\infty(m,i)$ as well): 
$$X_\infty(m,i)\subset \n P^m(\overline{\n F_q})$$

It turns out that the codimension of $X_\infty(m,i)$ in $\n P^m$ is much less than the quantity of these $\overline{H_{\be\al}}$. Further, it turns out that not only $\overline{H_{\be\al}}$, but also $H_{\be\al}$ are highly dependent. 
\medskip
In order to understand this phenomenon, we should consider not only the Carlitz module $\g C$ itself, but also its tensor powers and their twists. Namely, for any $n>0$ there exists $\g C^{\otimes n}$ --- the $n$-th tensor power of $\g C$. The formula for $L(\goth C^{\otimes n}_P,T)$ --- the L-series of its twist --- has a form similar to the one of $L(\goth C_P,T)$. We have the corresponding objects $\g M(\g a_0, \dots, \g a_m)(n)$, $L_m(n)$, $H_{\be\al}(n)$ etc. The value of $k$ for the case of any $n$ is $k=\lceil\frac{m+n}{q-1}\rceil$ (the analog of (3.1) for any $n$). 
\medskip
Since $L_m(n)\in\n Z[\g a_0, \dots, \g a_m,t,T]$, we can consider not only $$X_\infty(m,i)(n)\subset \n P^{m}(\overline{\n F_q})$$ which is the set of zeroes of $\overline{H_{\be\al}(n)}$ for $\be>k-i$, but its characteristic 0 analog (subscript $c0$ means "characteristic 0"): 
\medskip
$X_\infty(m,i)(n)_{c0}\subset \n P^{m}(\n C)$ which is the set of zeroes of $H_{\be\al}(n)$ for $\be>k-i$.
\medskip
Further, we can:
\medskip
First, to take $n=0$; the obtained matrix $\g M(\g a_0, \dots, \g a_m)(0)$ has entries in $\n Z[\g a_0, \dots, \g a_m]$, i.e. it does not contain $t$;
\medskip
Second, we can choose $q=2$. Although, as it was written earlier, for $q=2$ we have: $\g C^{\otimes n}$ have no forms (and hence in characteristic 2 the functions $L(\goth C^{\otimes n}_P,T)$ do not depend on $P$, up to finitely many Euler factors), but in characteristic 0 we have a non-trivial theory. 
\medskip
Until the end of Section 4, we give results for the case $q=2$. First of all, we have a non-expected fact: 
\medskip
{\bf Conjecture 4.1} (EGL], Conjecture 9). Let $q=2$. We have: the sets $X_\infty(m,i)(n)_{c0}$ do not depend on $n$, including $n=0$ (but they are different as schemes). 
\medskip
This means that $\forall \ n>0$ the polynomials $H_{\be\al}(n)$ belong to the radical of an ideal generated by some $H_{\be0}(0)$. See [EGL], Conjecture 9 for an exact, short, elementary self-contained statement of this result for $n=1$ (the most important case). 
\medskip
{\bf Remark 4.2.} Polynomials $H_{\be0}(0)$ are denoted in [EGL] by $\pm D(m,m-\be)$, and for $n=1$ the $H_{\be\al}$ of the present paper is denoted by $H_{\be, m-\al}(m)$ in [EGL], sorry  (all $H_{\be0}$ depend on $m$ and $n$, for $n=0$ the value of $\al$ in $H_{\be\al}(0)$ is always 0). The conjecture states that $\forall \ \al, \ \be,\ m, \ n$ there exists a power of $H_{\be\al}(m,n)$ (denoted by $\vk$) which belongs to the ideal generated by $H_{0,0}(m,0), H_{0,1}(m,0), \dots, H_{0,m-\be}(m,0)$. The conjecture is proved for some cases where $\vk=1$. Also, [EGL] contains a table for values of $\vk$ for $n=1$ and some small values of $\al, \ \be,\ m$. The values of $\vk$ given in this table look rather mysterious: we see few laws for them. Finally, [EGL] contains a numerical example of a representation of a $H_{\be\al}^2$ as a linear combination of $H_{0,*}(m,0)$ (the simplest example for the case $\vk=2$). 
\medskip
It was written above that matrices $\g M(\g a_0, \dots, \g a_m)$ are versions of a Schur matrix. It turns out that matrices $\g M(\g a_0, \dots, \g a_m)(0)$, that give rise to $H_{\be0}(0)=\pm D(m,m-\be)$, are exactly (up to a symmetry) Schur matrices for $\la_1=m, \ \la_2=m-1, \dots$ and $c_*$ from [FP], (1,5) are $\g a_*$. 
\medskip
The subject of [GLZ22], [GLZ23], [GLZ] is study of $X_\infty(m,i)(0)_{c0}$, case $q=2$. For simplicity, these varieties are denoted by $X(m,i)$ (do not confuse them with $X(m,i)$ of Section 3). Let us give the results of [GLZ22], [GLZ23], [GLZ] (some of them are conjectural).
\medskip
{\bf 4.3.} The set of irreducible components of $X(m,i)$ (denoted by $Irr(X(m,i))$ is described in terms of a set of rooted weighted binary trees and their disjoint unions (called forests) having $i$ nodes, up to some equivalence coming from a group action on the set of these forests. Weights (numbers attached to any node of a forest) are roots of 1 of order $2^\ga$ (here $\ga=0,\ 1,\ 2,\dots$) satisfying some relations. 
\medskip
More exactly, a description of $Irr(X(m,i))$ requires two consecutive constructions: the first construction (Section 5 of [GLZ22], [GLZ]) describes the so-called minimal components, and the second construction (Section 7 of [GLZ22], [GLZ]) describes how to get other components starting from the minimal ones. The second construction depends essentially on Conjecture 4.1 (this is astonishing, because apparently there is no relations between these subjects). 
\medskip
The result of [GLZ22], [GLZ] is not complete. We proved that a weighted forest gives us an irreducible component of $X(m,i)$; we conjecture that all irreducible components of $X(m,i)$ are obtained by the above two constructions. This conjecture is strongly supported by computer calculations. For example, the degree of $X(m,i)$ as a scheme is $$(m-1)(m-2)\dots(m-i)$$because $X(m,i)$ is a complete intersection (conjecturally; we cannot prove this fact) of $i$ hypersurfaces of degrees $m-1, \ m-2, \dots, \ m-i$. 
\medskip
Let $Y$ be an irreducible component of $X(m,i)$, $d(Y)$ its degree (as a subvariety of the ambient space $\n P^m$) and $\mu(Y)$ its multiplicity in $X(m,i)$. Therefore, $$\sum_{Y\in Irr(X(m,i))}d(Y)\mu(Y)=(m-1)(m-2)\dots(m-i).\eqno{(4.4)}$$
\medskip
We checked that (4.4) really holds for small $i$ and all $m$ if we sum over all $Y$ obtained by the first and second constructions, i.e. the description of $Irr(m,i)$ for these cases is exhaustive.
\medskip
{\bf 4.5.} All irreducible components of $X(m,i)$ are rational varieties. Their rational parametrizations are found. Hence, we can calculate their degrees $d(Y)$, fields of definitions (they are cyclotomic extensions of $\n Q$), actions of Galois groups on the set of forests, and partitions of $i$ naturally attached to any irreducible component $Y$ of $X(m,i)$.
\medskip
{\bf Problem 4.6.} There exists a simple formula for the multiplicity $\mu(Y)$ of any irreducible component $Y$ of $X(m,i)$. It is found experimentally, as a result of computer calculations. It should be proved. We have only a proof that if this formula is true for $Y$ obtained from some forests, then it is true for $Y$ obtained from some other forests.
\medskip
{\bf 5. Rank at infinity, case of any $q$.}
\medskip
Some results for this case are obtained in [GL16], Chapter 8. For $q>2$ we should distinguish the cases of finite characteristic and characteristic 0. 
\medskip
The principal obtained result is [GL16], Theorem 8.6. It states that 
$X(q,n,m,i)$ has at least two irreducible components: the first one, denoted by $X_r(q,m,i)$, does not depend on $n$, it has the codimension $i(q-1)$; the second one, denoted by $X_l(q,n,m,i)$, has the codimension $\le i(q-1)+(q-2)(n-1)$. Particularly, this means that polynomials $H_{\al\be}(n)$ are highly dependent. 
\medskip
Some examples for $q=3$, $n=1$, and small values of $m$ and $i$ are treated. Varieties $X(q,n,m,i)$ are interpreted in terms of determinantal varieties (see, for example, [FP] for their properties). 
\medskip
We should consider these results as an analog of the first theorems of the above theory for the case $q=2$, characteristic 0. Most likely $X(q,n,m,i)$ have much more irreducible components than the mentioned $X_r(q,m,i)$ and $X_l(q,n,m,i)$. An analog of Conjecture 4.1 should hold. Preliminary considerations show that to classify the irreducible components of $X(q,n,m,i)$ we should consider $q$-ary trees and forests, instead of binary ones for the case $q=2$. 
\medskip
{\bf 6. Possibilities of further research.} 
\medskip
First of all, we should prove Conjecture 4.1 and the fact that the description of irreducible components of $X(m,i)$ given in [GLZ22] is exhaustive. 
\medskip
Problems that are independent on these ones: 
\medskip
We do not know a description of singularities of irreducible components of $X(m,i)$ (they are highly singular) and their desingularizations. 
\medskip
There are tautological sheaves on irreducible components of $X(m,i)$. Their origin is the following. We can consider the matrix $\g M(\g a_0, \dots, \g a_m)$ as a linear transformation. There exist sheaves that correspond to kernels and images of its powers. Since irreducible components of $X(m,i)$ are rational varieties, these sheaves have a good description. We should find it. 
\medskip
Maybe the theory that we are going to develop will be an analog of the theory of Schubert subvarieties of flag varieties: Schubert subvarieties are also singular, we know a description of their singularities and tautological sheaves. 
\medskip
Further, we should get analogous results for $q>2$. An analog of Conjecture 4.1 for this case exists ([GL16], Theorem 8.6 and the below considerations), but its statement can be slighly different. 
\medskip
Finally, we want to pass from the order of zero of $L(\g C_P,T)$ at $T=\infty$ to the order of zero at $T=1$, in order to solve Problems 3.6.1 --- 3.6.4, and the principal Problems 1.4 and 3.4.
\medskip
Also, it is possible to develop the whole above theory for the case of any Drinfeld module instead of the Carlitz module. 

\medskip
\medskip
{\bf References}
\medskip

[EGL] S. Ehbauer, A. Grishkov, D. Logachev. Explicit statement of a conjecture on resultantal varieties. 2022.  https://arxiv.org/pdf/2209.04044.pdf
\medskip
[ELS] S. Ehbauer, D. Logachev, M. Sarraff de Nascimento. Some cases of a conjecture on L-functions of twisted Carlitz modules. Comm. Algebra 46 (2018), no. 5, 2130 -- 2145. https://arxiv.org/pdf/1707.04339.pdf
\medskip
[FP] Fulton, W., Pragacz, P. (1998). Schubert Varieties and Degeneracy Loci, Lecture Notes in Mathematics, 1689. Berlin: Springer-Verlag.
\medskip
[GL16] A. Grishkov, D. Logachev. Resultantal varieties related to zeroes of L-functions of Carlitz modules. Finite Fields and Their Applications. 2016, vol. 38, p. 116 -- 176. https://arxiv.org/pdf/1205.2900.pdf
\medskip
[GLZ22] A. Grishkov, D. Logachev, A. Zobnin. L-functions of Carlitz modules, resultantal varieties and rooted binary trees - I.  J. of Number Theory,  2022, vol. 238, p. 269 -- 312. 
\medskip
[GLZ23] Grishkov A., Logachev D., Zobnin A. L-functions of Carlitz modules, resultantal varieties and rooted binary trees - II. J. Algebra Appl. 22 (2023), no. 6, Paper No. 2350125
\medskip
[GLZ] is a complete version of the union of [GLZ22], [GLZ23]. https://arxiv.org/pdf/1607.06147.pdf 
\medskip
[PPVW] Park, Jennifer; Poonen, Bjorn; Voight, John; Wood, Melanie Matchett. A heuristic for boundedness of ranks of elliptic curves. J. Eur. Math. Soc. (JEMS) 21 (2019), no. 9, 2859 -- 2903
\medskip
\enddocument